\documentclass{article}
\usepackage{graphicx}
\usepackage{amsthm}
\usepackage{amsmath}
\usepackage{amsfonts}
\usepackage{amssymb}
\newtheorem{theorem}{Theorem}

\newtheorem{corollary}{Corollary}

\newtheorem*{definition}{Definition}

\newtheorem{proposition}{Proposition}

\begin{document}

\title{A note on the classification of naturally graded Lie algebras with linear
characteristic sequence}
\author{Jos\'{e} Mar\'{\i}a Ancochea Berm\'{u}dez\thanks{corresponding author: e-mail: Jose\_Ancochea@mat.ucm.es}
\and Rutwig Campoamor\thanks{e-mail: rutwig@nfssrv.mat.ucm.es}\\Departamento de Geometr\'{\i}a y Topolog\'{\i}a\\Fac. CC. Matem\'{a}ticas Univ. Complutense\\28040 Madrid ( Spain )}
\date{}
\maketitle

\begin{abstract}
For sufficiently high dimensions, the naturally graded nonsplit nilpotent Lie algebras with linear characteristic sequence are classified.
\end{abstract}

\section*{Introduction}

The determination by Vergne of the naturally graded filiform Lie algebras [6]
in relation with the study of irreducible components of the variety of
nilpotent Lie algebra laws pointed out the interest of graded structures for
this class of algebras. The introduction of the characteristic sequence of a
nilpotent Lie algebra [1] made it feasible to formulate the problem for
nilpotent Lie algebras of lower nilindexes.\newline Computational methods have
been applied in [5] to determine the algebras of characteristic sequence
$\left(  n-3,1,1,1\right)  $ for $n\geq8$. With few exceptions ( due to the
weakness of the conditions imposed by the graduation in low dimensions ) these
algebras belong to a couple of families. However, the weight of this procedure
relies heavily on the computations, which restricts the generality of the
results obtained. Following the philosophy of [6], i.e., using both extension
theory and the cohomology spaces with trivial coefficients, the problem can be
satisfactorily solved in arbitrary dimension and nilindexes, up to a reduced
number of pathological cases appearing in low dimensions. In this sense, all
the naturally graded Lie algebras of characteristic sequence $\left(
2m-1,q,1\right)  $ with $m\geq3$ and $1\leq q\leq2m-3$ were classified in
[2].\newline Finally. as a global classification of nilpotent Lie algebras
is a hopeless problem [4], it seems more reasonable to analyze the graded
structures and then use deformation theory to obtain concrete isomorphism
classes satisfying certain properties. This has for example been done to
obtain nilradicals of complete and rigid Lie algebras [3].

\section{Preliminaries and notations}

Whenever we speak about Lie algebras in this work, we refer to finite
dimensional complex nonsplit Lie algebras.\bigskip

Let $\frak{g}$ be an $n$-dimensional nilpotent Lie algebra. We denote by
$C^{i}\left(  \frak{g}\right)  $ the ideals of the central descending sequence
of $\frak{g}$. Then the $\mathbb{Z}$-graded Lie algebra associated to this
filtration is defined by
\[
\frak{gr}\left(  \frak{g}\right)  =\sum_{i\in\mathbb{Z}} \frac{C^{i}\frak{g}%
}{C^{i+1}\frak{g}}
\]

\begin{definition}
A nilpotent Lie algebra $\frak{g}$ is called naturally graded is it is
isomorphic to $\frak{gr}\left(  \frak{g}\right)  $.
\end{definition}

In [1] the following invariant is presented :\newline 

For any nonzero vector $X\in\frak{g}_{n}-C^{1}\frak{g}_{n}$ let $c\left(
X\right)  $ be the ordered sequence of a similitude invariant for the
nilpotent operator $ad_{\mu}\left(  X\right)  \dot{,}$ i.e., the ordered
sequence of dimensions of Jordan blocks of this operator. The set of these
sequences is ordered lexicographically.

\begin{definition}
The characteristic sequence of $\frak{g}_{n}$ is an isomorphism invariant
$c\left(  \frak{g}_{n}\right)  $ defined by
\[
c\left(  \frak{g}_{n}\right)  =\max_{X\in\frak{g}_{n}-C^{1}\frak{g}_{n}%
}\left\{  c\left(  X\right)  \right\}
\]
A nonzero vector $X\in\frak{g}_{n}-C^{1}\frak{g}_{n}$ for which $c\left(
X\right)  =c\left(  \frak{g}_{n}\right)  $ is called characteristic
vector.\newline A characteristic sequence is called linear if it is of the
form $\left(  q,1,..,1\right)  $ for $q\geq2$.
\end{definition}

A Lie algebra $\frak{g}$ is called filiform if $c\left(  \frak{g}\right)
=\left(  n-1,1\right)  $, where $n$ is the dimension. In [6], Vergne proves
that every naturally graded filiform Lie algebra is isomorphic to one of the
algebras $L_{n}$ and $Q_{n}$ defined by

\begin{enumerate}
\item $L_{n}$ $\left(  n\geq3\right)  :$
\[
\lbrack X_{1},X_{i}]=X_{i+1},\;1\leq i\leq n
\]
over the basis $\left\{  X_{1},..,X_{n+1}\right\}  $.

\item $Q_{2m-1}$ $\left(  m\geq3\right)  :$%
\begin{align*}
\lbrack X_{1},X_{i}] &  =X_{i+1},\;1\leq i\leq2m-1\\
\lbrack X_{j},X_{2m+1-j}] &  =\left(  -1\right)  ^{j}X_{2m},\;1\leq j\leq m
\end{align*}
over the basis $\left\{  X_{1},..,X_{2m}\right\}  $.
\end{enumerate}

It will be convenient to use the so called contragradient representation of a
Lie algebra $\frak{g}$. Let $n=dim\left(  \frak{g}\right)  $ and $\left\{
X_{1},..,X_{n}\right\}  $ be a basis. If $C_{i,j}^{k}$ are the structure
constants of the algebra law $\mu$, we can define, over the dual basis
$\left\{  \omega_{1},..,\omega_{n}\right\}  $, the differential
\[
d_{\mu}\omega_{i}\left(  X_{j},X_{k}\right)  =-C_{jk}^{i}
\]
Then the Lie algebra is rewritten as
\[
d\omega_{i}=-C_{jk}^{i}\omega_{j}\wedge\omega_{k}\; 1\leq i,j,k\leq n,
\]
The Jacobi condition is equivalent to $d^{2}\omega_{i}=0$ for all $i$. Thus
the problem of determining the structure constants is translated to the
calculation of a system of closed differential forms. This approach is more
useful in practice, as it allows to recognize both extensions and central
quotients faster. This also explains why for the first dimensions admitting a
certain graduation type and a linear characteristic sequence, the conditions
on the structure constants are too weak.

The two central results to classify the naturally graded Lie algebras in
dimension $n$ and linear characteristic sequence $\left(  n-t,1,..,1\right)  $
for $2\leq t\leq\left[  \frac{n}{2}\right]  $ are

\begin{proposition}
Let $\frak{g}$ be a split naturally graded Lie algebra of nilindex $p$ and
linear characteristic sequence. Then $\frak{g}$ does not admit a nonsplit
naturally graded central extension of degree one preserving both the nilindex
and the linearity of $c\left(  \frak{g}\right)  $.
\end{proposition}

\begin{proposition}
Let $\frak{g}$ an $n$-dimensional naturally graded Lie algebra with nilindex
$p$ and linear characteristic sequence. Then there exists an $\left(
n-1\right)  $-dimensional naturally graded Lie algebra $\frak{g}^{\prime}$
with the same nilindex and linear characteristic sequence such that the
following sequence is exact :
\[
0\rightarrow\mathbb{C}\rightarrow\frak{g}\rightarrow\frak{g}^{\prime
}\rightarrow0
\]
\end{proposition}

The proof of the first assertion is trivial, as any nonsplit graded extension
preserving the nilindex and the linearity of the characteristic sequence must
necessarily involucrate the characteristic vectors, which implies the
generation of a Jordan block of dimension two for the adjoint operator of such
a vector. As to the second proposition, for sufficient large dimensions it is
not difficult to see that the conditions on the coefficients of the
contragradient representation do not depend any more on the dimension of the
algebra, but only on the graduation type [2].\newline Starting from these
auxiliary results and proceeding like in [2] we obtain :

\begin{theorem}
For sufficiently high dimension, any naturally graded Lie algebra $\frak{g}$
of characteristic sequence $\left(  dim\frak{g}-k,1,.\overset{\left(
k\right)  }{}.,1\right)  $ with $2\leq k\leq\left[  \frac{\dim\frak{g}}%
{2}\right]  $ is isomorphic to one of the following models :

\begin{enumerate}
\item $L_{n}\left(  t_{1},..,t_{p}\right)  \,;\;1\leq t_{1}<..<t_{p}%
\leq\left[  \frac{n-1}{2}\right]  ,\;n\geq 7$%
\[%
\begin{tabular}
[c]{l}%
$d\omega_{1}=d\omega_{2}=0$\\
$d\omega_{j}=\omega_{1}\wedge\omega_{j-.1},\;3\leq j\leq n$\\
$d\omega_{n+i}=$\ $\displaystyle\sum_{j=2}^{t_{i}+1}\left(  -1\right)
^{j}\omega_{j}\wedge\omega_{2t_{i}+3-j},\;1\leq i\leq p$%
\end{tabular}
\]

\item $Q_{2m-1}\left(  t_{1},..,.t_{p}\right)  \;;\;1\leq t_{1}<..<t_{p}%
\leq\left[  \frac{2m-1}{2}\right] ,\;m\geq 3 $%
\[%
\begin{tabular}
[c]{l}%
$d\omega_{1}=d\omega_{2}=0$\\
$d\omega_{j}=\omega_{1}\wedge\omega_{j-1},\;3\leq j\leq2m-1$\\
$d\omega_{2m}=\omega_{1}\wedge\omega_{2m-1}+\displaystyle\sum_{j=2}^{\left[
\frac{2m+1}{2}\right]  }\left(  -1\right)  ^{j}\;\omega_{j}\wedge
\omega_{2m+1-j}$\\
$d\omega_{2m+1+j}=\displaystyle\sum_{i=2}^{t_{j}}\left(  -1\right)  ^{j}%
\omega_{i}\wedge\omega_{2t_{j}+3-i},\;1\leq j\leq p$%
\end{tabular}
\]

\item $D_{2m}\left(  t_{1},..,t_{p}\right)  \;;\;1\leq t_{1}<..<t_{p}\leq m-3 ,\;m\geq 4$%
\[%
\begin{tabular}
[c]{l}%
$d\omega_{1}=d\omega_{2}=0$\\
$d\omega_{j}=\omega_{1}\wedge\omega_{j-1},\;3\leq j\leq2m-3$\\
$d\omega_{2m-2}=\omega_{1}\wedge\omega_{2m-3}+\displaystyle\sum_{j=2}%
^{m-1}\left(  -1\right)  ^{j}\;\omega_{j}\wedge\omega_{2m-3-j}$\\
$d\omega_{2m-1}=\omega_{1}\wedge\omega_{2m-2}+\displaystyle\sum_{j=2}%
^{m-1}\left(  -1\right)  ^{j}\left(  m-j\right)  \,\ \omega_{j}\wedge
\omega_{2m-2-j}$\\
$+\left(  2-m\right)  \,\omega_{2}\wedge\omega_{2m}$\\
$d\omega_{2m}=\displaystyle\sum_{j=2}^{m-1}\left(  -1\right)  ^{j}\;\omega
_{j}\wedge\omega_{2m-3-j}$\\
$d\omega_{2m+i}=\displaystyle\sum_{j=2}^{t_{i}-1}\left(  -1\right)
^{j}\;\omega_{j}\wedge\omega_{2t_{i}+3-j},\;1\leq i\leq p$%
\end{tabular}
\]

\item $E_{2m+1}\left(  t_{1},..,t_{p}\right)  \;;\;1\leq t_{1}<..<t_{p}\leq
m-3 ,\;m\geq 4$%
\[%
\begin{tabular}
[c]{l}%
$d\omega_{1}=d\omega_{2}=0$\\
$d\omega_{j}=\omega_{1}\wedge\omega_{j-1},\;3\leq j\leq2m-3$\\
$d\omega_{2m-2}=\omega_{1}\wedge\omega_{2m-3}+\displaystyle\sum_{j=2}%
^{m-1}\left(  -1\right)  ^{j}\;\omega_{j}\wedge\omega_{2m-1-j}$\\
$d\omega_{2m-1}=\omega_{1}\wedge\omega_{2m-2}+\displaystyle\sum_{j=2}%
^{m-1}\left(  -1\right)  ^{j}\left(  m-j\right)  \,\ \omega_{j}\wedge
\omega_{2m-j}$\\
$+\left(  m-2\right)  \,\omega_{2}\wedge\omega_{2m+1}$\\
$d\omega_{2m}=\omega_{1}\wedge\omega_{2m-1}+\displaystyle\sum_{j=3}^{m}%
\frac{\left(  -1\right)  ^{j}\left(  j-2\right)  \left(  2m-1-j\right)  }%
{2}\;\omega_{j}\wedge\omega_{2m+1-j}$\\
$+\left(  m-2\right)  \,\omega_{3}\wedge\omega_{2m+1}$\\
$d\omega_{2m+1}=\displaystyle\sum_{j=2}^{\left[  \frac{2m+1}{2}\right]
-1}\left(  -1\right)  ^{j}\;\omega_{j}\wedge\omega_{2m-1-j}$\\
$d\omega_{2m+1+i}=\displaystyle\sum_{j=2}^{t_{i}+1}\left(  -1\right)
^{j}\omega_{j}\wedge\omega_{2t_{i}+3-j},\;1\leq i\leq p$%
\end{tabular}
\]
\end{enumerate}
\end{theorem}

The result can be resumed in the following table. The fourth column indicates
the restrictions on $n$.
\[%
\begin{tabular}
[c]{|c|c|c|c|}\hline
$\frak{g}$ & $\dim\;\frak{g}$ & $ch\left(  \frak{g}\right)  $ & \\\hline
$L_{n}\left(  t_{1},..,t_{p}\right)  $ & $n+p$ & $\left(  n-1,1,..^{\left(
p+1\right)  },..,1\right)  $ & \\\hline
$Q_{n}\left(  t_{1},..,t_{p}\right)  $ & $2m+p$ & $\left(  2m-1,1,..^{\left(
p+1\right)  },..,1\right)  $ & $n=2m-1$\\\hline
$D_{n}\left(  t_{1},..,t_{p}\right)  $ & $2m+p$ & $\left(  2m-2,1,..^{\left(
p+2\right)  },..,1\right)  $ & $n=2m$\\\hline
$E_{n}\left(  t_{1},..,t_{p}\right)  $ & $2m+1+p$ & $\left(
2m-1,1,..^{\left(  p+2\right)  },..,1\right)  $ & $n=2m+1$\\\hline
\end{tabular}
\]

\begin{corollary}
For any dimension $n\geq7$ the number of naturally graded Lie algebras with
linear characteristic sequence is finite.
\end{corollary}

\end{document}